%

\magnification=1200
\nopagenumbers
\hsize=6truein
\hoffset=.2truein
\vsize=8truein
\baselineskip=14truept

\font\claim=cmcsc10 at 12truept
\font\smallit=cmti10 at 10truept

\def\today{\ifcase\month\or
    January\or February\or March\or April\or May\or June\or
        July\or August\or September\or October\or November\or December\fi
            \space\number\day, \number\year}

\def\hcm#1{\hskip #1truecm}
\def\hpt#1{\hskip #1truept}

    \newcount\referno
    \referno=1
\def\ref#1#2#3#4{
\item{[\the\referno]} {{\claim #1} \hpt 9 {\sl #2} \hpt 9 {\bf #3} \hpt
7 {#4}} \advance \referno by 1\vskip 2truept}

\def\reference{\bigskip{\claim References}\par\vskip 20truept}

\pageno=1
\headline={
\ifnum\pageno=1\hfill \else
\hfill  \quad
       \smallit \folio
\fi}

    \newcount\ssno
    \ssno=1

    \newcount\defno
    \defno=1
\def\definition{ \medskip \noindent
{\claim Definition \the\ssno .\the\defno \hcm 1}   \advance \defno by 1 }

    \newcount\exmno
    \exmno=1
\def\example{ \medskip \noindent
{\claim Example \the\ssno .\the\exmno \hcm 1}   \advance \exmno by 1 }

\def\Definition#1{\medskip\noindent{\claim Definition
     \the\ssno .\the\defno \hcm 1}   \advance \defno by 1 }

    \newcount\lemmano
    \lemmano=1
\def\lemma{ \medskip \noindent {\claim Lemma
\the\ssno .\the\lemmano\hcm 1}
 \advance \lemmano by 1}

\def\Lemma#1{\medskip \noindent{{\claim Lemma}
\the\ssno .\the\lemmano\hpt 9 #1 \hcm 1 } \advance \lemmano by 1}

    \newcount\thmno
    \thmno=1
\def\theorem{ \medskip \noindent {\claim Theorem
\the\ssno .\the\thmno\hcm 1}
         \advance \thmno by 1}
\def\Theorem#1{\medskip \noindent {\claim Theorem
\the\ssno .\the\thmno \hpt 5 #1 \hcm 1 } \advance \thmno by 1}

    \newcount\corno
    \corno=1
\def\corollary{ \medskip \noindent {\claim Corollary
\the\ssno .\the\corno\hcm 1}
         \advance \corno by 1}
\def\Corollary#1{\medskip \noindent {\claim Corollary
\the\ssno .\the\corno \hpt 5 #1 \hcm 1 } \advance \corno by 1}

\def\proof{\medskip \medskip \noindent{\bf Proof \hpt 9}}

\def\st#1{\vskip .4truein
       \noindent{ \S   \the\ssno \hskip
       1truecm {\claim #1} }  \vskip .1truein \noindent }

\def\endsection{\advance \ssno by 1
        \defno=1
        \lemmano=1
        \thmno=1
        \corno=1}

\chardef\lprime="12

\def\rtr{{\hbox{\rm $|$\kern-2.6pt\lprime}}}

\def\forces{\thinspace \hbox{{\vrule height6.9pt depth0.0pt
width0.2pt}\kern0.4pt$\vdash$ }}

\def\square{ \vbox{\hrule\hbox{\vrule\kern3pt
                    \vbox{\kern3pt\kern3pt}\kern3pt\vrule}\hrule}}


\ssno=0

\st{Introduction}

In [3], Foreman, Magidor and Shelah introduced a forcing axiom, called
Martin's Maximum, which says that if $P$ is a forcing notion preserving
stationary subsets of $\omega_1$, and $\langle D_\alpha \;|\;\alpha <
\omega_1\rangle$ is a sequence of dense subsets, then there exists a
filter $G \subseteq P$ meeting every $D_\alpha.$ They showed that if the
existence of a supercompact cardinal is consistent with the set theory,
then Martin's Maximum is also consistent with the set theory. They also
showed many consequences of Martin's Maximum. For example, Martin's
Maximum implies that the Singular Cardinal Hypothesis holds, the
nonstationary ideal on $\omega_1$ is saturated, and every stationary
set on every countable space reflects.

In this paper, we continue the study of Martin's Maximum initiated in
[3]. We define {\it projective} stationary sets and
prove that Martin's Maximum implies that every projective
stationary set contains an increasing continuous $\in$--chain of length
$\omega_1$. We will also show that the consequences of Martin's Maximum
in [3] follow from the assumption that every projective stationary set
contains an increasing continuous $\in$--chain of length $\omega_1$.  We
will call
this assumption the Projective Stationary Reflection principle. It
turns out that this principle is
in fact equivalent to the Strong Reflection Principle of Todorcevic [9]
(see next section for details).
It differs from the known stationary reflection principles in two
respects. First, it is not applicable to every stationary set, but the
projective stationary sets. Basically,  a projective stationary set is a
stationary set which projects to every stationary subset of $\omega_1$.
Hence there are many stationary sets which are not projective
stationary. Secondly, it claims that every projective stationary set
reflects to a closed unbounded subset, not just a stationary subset
, of some countable space $[X]^\omega$  with $X$ has cardinality of
$\aleph_1$.

It is a strong reflection principle because it implies the known
stationary reflection principles, not the other way around. A typical
stationary reflection principle is the following reflection principle
from [3]:

\vskip 5truept
\noindent{(RP)}\hskip 1truecm {\sl If $S \subseteq [H_\kappa]^{\omega}$
is stationary, then there exists an $X$ of size $\aleph_1$ such that
$\omega_1 \subseteq X$ and $S\cap [X]^{\omega}$ is stationary in
$[X]^{\omega}$.}
\vskip 5truept
This stationary reflection  principle has many interesting consequences.
For example, RP implies that every $\omega_1$--stationary preserving
partially ordered set is semiproper;
RP implies that the nonstationary ideal on $\omega_1$ is
presaturated (see later for definitions), theorems of Foreman,
Magidor and Shelah [3]; RP implies that $2^{\aleph_0} \leq \aleph_2$, a
theorem
of Todorcevic [8],
etc. (See [1,2,3,8,9,10] for more results on this
direction.)

However, it is also proved in [3] that when a supercompact cardinal is
Levy
collapsed to $\omega_2$, in the resulting model, every stationary set
reflects. But in such a model, the strong reflection principle does not
hold.

In section 1, we study the Projective Stationary Reflection Principle.
We will prove that it follows from Martin's Maximum and it is
equivalent to the Strong Reflection Principle of Todorcevic.

In section 2, we present some natural examples of projective stationary
sets.

In section 3, we study some special projective stationary sets, which
form a filter. It turns out that this filter is closely connected to the
saturation of the nonstationary ideal on $\omega_1$. Namely, the
saturation of the nonstationary ideal on $\omega_1$ itself is a kind of
reflection.

We refer to [4,5] for all terms used but not explicitly defined in the text.

\endsection

\st{Projective Stationary Sets }

 Let $\kappa$ be a regular cardinal $\geq \omega_2$. Let
$H_\kappa$ be the set of all sets hereditarily of size less than
$\kappa$. We assume that $H_\kappa$ is endowed with a well ordering.
We use $[H_\kappa]^{\omega}$ to denote the set of all
countable subsets of $H_\kappa$. Notice that $[H_\kappa]^\omega
\subseteq H_\kappa.$ A subset $C \subseteq
[H_\kappa]^{\omega}$ is {\bf closed} and {\bf unbounded} (a {\bf club},
in short)
if there is a function $f :[H_\kappa]^{<\omega} \to H_\kappa$ such that
for every countable $x \in [H_\kappa]^{\omega}$, $x \in C$ if and only
if $x$ is closed under $f$, i.e., if $e \in [x]^{<\omega}$, then $f(e)
\in x$, where $[A]^{<\omega}$ denotes the set of all finite subsets of
$A$. A subset $S \subseteq [H_\kappa]^{\omega}$ is {\bf stationary} if
for every club $C \subseteq [H_\kappa]^{\omega}$ the intersection $C\cap
S$ is not empty. A basic fact that is used frequently is the normality 
of the club
filter: if
$S$ is stationary, $f:S \to H_\kappa$ is a choice function, i.e., $f(x)
\in x$ for all $x \in S,$ then there is a stationary subset on which $f$
is a constant. In particular, the intersection of countably many clubs
contains a club.

We use $N \prec M$ to denote that $N$ is an elementary submodel of $M$
as usual. Also if $f : X \to Y$ and $A \subseteq X$, we use $f''A$ to
denote the set $\{f(a)\;|\;a \in A\}.$

A sequence $\langle N_\alpha\;|\;\alpha < \theta\rangle$ of countable
submodels of $H_\kappa$ is an {\bf increasing continuous
$\in$--chain} of length
$\theta$ if for all $\alpha < \theta,\;(N_\alpha\prec H_\kappa)$
and for all $\alpha < \beta < \theta,\;(N_\alpha \in N_\beta),$
and if $\alpha < \theta$ is a limit ordinal then $N_\alpha =
\bigcup_{\beta<\alpha}N_\beta.$

A sequence $\langle N_\alpha\;|\;\alpha < \theta\rangle$ of countable
submodels of $H_\kappa$ is a {\bf strongly increasing continuous
$\in$--chain} of length
$\theta$ if for all $\alpha < \theta,\;(N_\alpha\prec H_\kappa)$
and for all $\alpha < \beta < \theta,\;(N_\alpha \in N_\beta),$
and if $\alpha < \theta$ is a limit ordinal then $N_\alpha =
\bigcup_{\beta<\alpha}N_\beta$ and $\alpha+1 < \theta$ implies that
$\langle N_\beta \;|\;\beta < \alpha \rangle \in N_{\alpha+1}.$

A sequence $\langle N_\alpha\;|\;\alpha < \theta\rangle$ of countable
submodels of $H_\kappa$ is an {\bf increasing continuous
$\subseteq$--chain} of length
$\theta$ if for all $\alpha < \theta,\;(\alpha \subseteq N_\alpha\prec
H_\kappa)$
and for all $\alpha < \beta < \theta,\;(N_\alpha \subseteq N_\beta),$
and if $\alpha < \theta$ is a limit ordinal then $N_\alpha =
\bigcup_{\beta<\alpha}N_\beta.$

We will frequently use a simple fact without mentioning it. This fact
says that if $N \prec H_\kappa$ and $x \in N$ is countable then $x
\subseteq N.$
It follows that every increasing continuous $\in$--chain is an
increasing continuous $\subseteq$--chain.

Note that for closed unbounded many $X \in [H_\kappa]^\omega,$
$X \cap \omega_1$ is an ordinal.

\definition The {\bf projection} of a set $S \subseteq [H_\kappa]^\omega$
is the set Proj$(S)=\{X \cap \omega_1 : X \in S\}.$
We say that $S \subseteq [H_\kappa]^{\omega}$
is {\bf
projective stationary} if for every 
club $C \subseteq [H_\kappa]^{\omega}$, Proj$(S \cap C)$ contains a club
in $\omega_1.$

\vskip 5truept
It follows from the definition that every projective stationary set is
stationary. But the converse is not true. Also, there are disjoint
projective stationary sets. We will see some natural examples later.

We are interested in projective stationary sets primarily because we
are interested in the following
Projective Stationary Reflection principle:

\vskip 5truept
\noindent{\bf Projective Stationary Reflection:} {\sl For every
regular
cardinal $\kappa \geq \omega_2$, if $S\subseteq [H_\kappa]^{\omega}$ is
a projective stationary
set, then there exists an increasing continuous $\in$--chain $\langle
N_\alpha\;|\;\alpha < \omega_1\rangle$ of countable elementary submodels
of $H_\kappa$ of length $\omega_1$ such that $N_\alpha \in S$ for all
$\alpha < \omega_1.$ }
\vskip 5truept
We now show that the Projective Stationary Reflection Principle follows
from the Martin's Maximum.
\vskip 5truept
Let us recall Martin's Maximum:

A $\omega_1$--stationary preserving partially ordered set is a partially
ordered set forcing over the ground model with it every stationary
subset of $\omega_1$ in the ground model remains to be stationary in
the generic extension. Martin's Maximum is the following statement:
\medskip
{\sl If $P$ is a $\omega_1$--stationary preserving partially
ordered set, if
$\{D_\alpha \;|\;\alpha < \omega_1\}$ is a sequence of dense subsets of
$P$ of size $\aleph_1$, then there exists a filter $G \subseteq P$ which
meets every $D_\alpha$.}
\medskip
It is proved in [3] that if the existence of a supercompact cardinal
is consistent with the axioms of set theory, then the Martin's
Maximum is also consistent with the axioms of set theory.
It is shown also in [3] that the Martin's Maximum has many
interesting consequences.

To show that the Martin's Maximum implies the Projective
Stationary Reflection Principle,
we need a lemma saying that every stationary subset $S \subseteq
[H_\kappa]^\omega$ contains strongly increasing continuous $\in$--chains
of any countable length.
This is a generalization of a well known fact that every stationary
subset of $\omega_1$ contains closed subsets  of any countable order
type (see [4], Exercise 7.12, Page 60).

\lemma If $S \subseteq \omega_1$ is stationary, then for any countable
ordinal $\alpha$, for any countable ordinal $\gamma$, there exists an
$x \subseteq S-\gamma$ such that $x$ is closed and has order type
$\alpha+1$.
\medskip
The proof of this lemma is by induction on $\alpha.$ We will use this
lemma to prove the following:

\lemma If $S \subseteq [H_\kappa]^\omega$ is stationary, then for every
countable
ordinal $\alpha$, there exists a strongly increasing continuous
$\in$--chain
$\langle N_\gamma\;|\;\gamma \leq \alpha\rangle$ of length $\alpha+1$
such that $N_\gamma \in S$ for all $\gamma \leq \alpha$.

\proof
Let $S \subseteq [H_\kappa]^\omega$ be stationary. Let $\alpha <
\omega_1$ be fixed. We will show that in some generic extension by a
$\sigma$--closed forcing the conclusion of the lemma holds for this $S$
and $\alpha$. That will be enough.

Consider the following forcing $P$.

A condition $p$ is a strongly increasing continuous $\in$--chain of
length $\theta+1$ for some $\theta < \omega_1$ such that $p(\beta) \prec
H_\kappa$ is countable for all $\beta < \theta+1.$

The ordering is by extension.
\smallskip
\noindent{\bf Claim}\hskip 1truecm $P$ is  $\sigma$--closed.
\smallskip
To see this, let $\langle p_n\;|\;n<\omega\rangle$ be a decreasing
sequence of conditions. Let $\theta_n$ be the largest countable ordinal
in the domain of $p_n$. Let $\theta = \bigcup_{n < \omega}\theta_n$ and
let $N = \bigcup_{n<\omega}p_n(\theta_n).$ Then $N \prec H_\kappa$ is a
countable elementary submodel. Define $q(\theta) = N$ and $q(\beta) =
p_n(\beta)$ for $\beta < \theta$ with $n$ the least such that $\beta
\leq \theta_n.$ $q$ is then a condition extending all the $p_n$.

\medskip
Hence $P$ is a $\sigma$--closed forcing.

Let $G\subseteq P$ be a generic
filter over $V$.  In $V[G]$, let $f = \bigcup G.$ Then $f$ is an
increasing continuous $\in$--chain of length $\omega_1$ which enumerates
a closed unbounded subset of $([H_\kappa]^\omega)^V$. Since $S$ remains
to be stationary in the generic extension, the following set is a
stationary subset of $\omega_1$:
$$T = \{\beta < \omega_1\;|\;f(\beta) \in S\}.$$
By the quoted lemma above, let $x \subseteq T$ be a closed set of order
type $\alpha+1$. Let $x = \{\gamma_\beta\;|\;\beta \leq \alpha\}$ be the
canonical enumeration. Let $N_\beta = f(\gamma_\beta)$ for $\beta \leq
\alpha.$ Then$\langle N_\beta \;|\;\beta \leq \alpha\rangle \in V$ is
the desired strongly increasing continuous $\in$--chain of length
$\alpha+1$.

This proves the lemma.

\hfill\square

\theorem Assume the Martin's Maximum. Then the Projective Stationary
Reflection principle holds. In fact, if $S \subseteq [H_\kappa]^\omega$
is projective stationary, then there exists a strongly increasing
continuous $\in$--chain of length $\omega_1$ through $S$.

\proof
Let $S \subseteq [H_\kappa]^\omega$ be projective stationary.

The idea is to shoot a strongly increasing continuous $\in$--chain
of length $\omega_1$ through $S$.  Thus, a
condition $p = \langle N_\alpha\;|\;\alpha \leq \theta\rangle$ is a
strong increasing continuous $\in$--chain of length $\theta+1 $ for some
$\theta < \omega_1$ such that for all
$\alpha \leq\theta\;(N_\alpha \in S).$    The ordering is by extension.

Let $P$ be the set of all conditions. We are going to show that forcing
with $P$ preserves stationary subsets of $\omega_1$.

First, we need to check that the following sets are dense.

For $\alpha < \omega_1$, let
$D_\alpha = \{ p \in P\;|\;\alpha \in dom(p)\}.$
For $x \in H_\kappa$, let $D_x = \{p \in P\;|\;\exists\,\alpha \in
dom(p)\;(x \in p(\alpha))\}.$

\medskip
We would like to show that each $D_\alpha$, $D_x$ is dense.
\medskip
To see this, let $\alpha < \omega_1$ and $x \in H_\kappa.$ Let $p =
\langle N_\beta\;|\;\beta \leq \theta\rangle$ be a condition.

Then $S_1 = \{ N \in S\;|\;\{p, x \} \subseteq N\}$ is
stationary. Applying the lemma above, let $\langle M_\gamma\;|\; \gamma
\leq \alpha \rangle$ be a strongly increasing continuous $\in$--chain
from
$S_1$. We then define $q(\beta) = M_\beta$ for all $\beta \leq \alpha$
such that $\beta > \theta$ and define $q(\beta) = p(\beta)$ for all
$\beta \leq \theta.$ It follows that $q \leq p$ and $q \in D_\alpha \cap
D_x.$

\medskip
We can now show that the forcing preserves stationary subsets of
$\omega_1$.

Let $T \subseteq \omega_1$ be stationary. Let $\dot C$ be a name for a
club subset of $\omega_1$. We would like to show that $\forces \check T
\cap\dot C \not=\emptyset.$

Let $p \in P$ be a condition.

Let $S_T = \{N\in S\;|\;N\cap \omega_1 \in T\}.$ Then $S_T$ is
stationary.

Let $\lambda \geq (2^{2^{|P|}})^+$ be a regular cardinal. Consider the
structure
$${\cal H} = \langle H_\lambda, \in,\triangle, P, \{\dot C\}, S_T,
\cdots\rangle.$$

Let $N \prec {\cal H}$ be countable such that $N\cap H_\kappa \in S_T$,
and  $\{\dot C, S, T, S_T, p \}\subseteq N.$ Let $\delta = N\cap
\omega_1.$ Let $\langle D_n\;|\;n < \omega\rangle$ be a list of all
dense subsets of $P$ which are in $N$. Let $p_0 = p.$ Inductively, pick
$p_{n+1} \in D_n\cap N$ so that $p_{n+1} \leq p_n.$ Let $\theta_n$ be
the
largest countable ordinal in $dom(p_n)$. By elementarity and a density
argument, we have $\delta = \bigcup_{n < \omega}\theta_n$ and $N\cap
H_\kappa = \bigcup_{n<\omega}p_n(\theta_n).$ Therefore,
$$q = \bigcup_{n<\omega}p_n \cup\{(\delta, N\cap H_\kappa)\}$$
is a condition stronger than all $p_n$. Since $\dot C \in N$, for each
$\beta < \delta$, there is a name $\dot \gamma \in N$ such that
$$\forces \dot \gamma \in \dot C\;\&\;\check \beta < \dot \gamma.$$
Each such name corresponding to a dense subset in $N$. Hence, $ q
\forces \dot\gamma \in \check \delta.$ It follows that $q\forces ``\dot
C \cap \check \delta$ is unbounded in $\check \delta.$''. Thus, $q
\forces \check \delta \in \dot C$. Therefore, $q \forces \dot C \cap
\check T \not= \emptyset.$

It follows then that forcing with $P$ preserves stationary subsets of
$\omega_1.$

\medskip
Now let $G \subseteq P$ be a filter meeting all the dense subsets
$D_\alpha$ for $\alpha < \omega_1$ defined above. Let
$$\langle N_\alpha\;|\;\alpha < \omega_1\rangle = \bigcup G.$$
Then this is a strongly increasing continuous $\in$--chain of length
$\omega_1$ with each $N_\alpha \in S.$

This finishes the proof.

\hfill\square
\medskip From the proof above, we have seen that for a given stationary 
$S\subseteq [H_\kappa]^\omega$, there is a natural forcing notion $P_S$
associated with $S$ to shoot an increasing continuous $\in$--chain of
length $\omega_1$. Then $S$ is projective stationary if and only if this
forcing notion $P_S$ preserves stationary sets of $\omega_1$.
\medskip

After the work of Foreman, Magidor and Shelah [3], Todorcevic in a
circulated hand written note [9] in September 1987 formulated the
following Strong Reflection Principle (SRP) (See also [1], page 57--60):
 \medskip
\noindent{\bf Strong Reflection Principle:} {\sl For every $\kappa$,
every $S \subseteq [\kappa]^{\omega}$ and
for every regular $\theta > \kappa$ there is an increasing continuous
$\in$--chain
$\{N_\alpha\;|\;\alpha < \omega_1\}$ of countable elementary models of
$H_\theta$ (with $N_0$ containing a prescribed element of $H_\theta$)
such that for all $\alpha < \omega_1$, $N_\alpha \cap \kappa \in S$ if
and only if there exists a countable elementary submodel $M$ of
$H_\theta$ such that $N_\alpha \subseteq M$, $M\cap \omega_1 =
N_\alpha\cap \omega_1$, and $M\cap \kappa \in S.$}
\medskip
In the following, we show that the Strong Reflection Principle and the
Projective Stationary Reflection principle are equivalent.

\theorem The Strong Reflection Principle holds if and only if the
Projective Stationary Reflection principle holds.

\proof
First we show that the Strong Reflection Principle implies the
Projective Reflection principle.

Let $S \subseteq [H_\kappa]^{\omega}$ be a projective stationary set.
Let $\{N_\alpha \;|\;\alpha < \omega_1\}$ be a continuous $\in$--chain
of countable elementary submodels of $H_\theta$ with $\theta > \kappa$,
given by SRP for $S$.

We would
like to show that $\{\alpha < \omega\;|\; N_\alpha\cap H_\kappa\in S\} $
contains a club in $\omega_1$.

Assume not. Let $T = \{\alpha < \omega_1\;|\;N_\alpha \cap H_\kappa
\not\in S \hbox{\rm \  and }N_\alpha\cap\omega_1 = \alpha \}$. Then $T$
is stationary.

 Define $D$ to be the following set
$$D =\{N\in [H_\theta]^{\omega}\;|\;H_\kappa \in N\;\&\;\forall \beta
\in N\cap \omega_1\;N_\beta \in N\}.$$
By normality, $D$ contains a club on $[H_\theta]^{\omega}$.

Since $S$ is projective stationary, we can have some $N \in D$
such that $N\cap \omega_1 \in T$, $N\cap H_\kappa \in S$ and $N$ is an
elementary submodel
of $H_\theta.$ Let $\alpha = N\cap \omega_1.$ Since if $\beta < \alpha$
then $N_\beta \subseteq N$, $N_\alpha \cap \omega_1 = \alpha = N \cap
\omega_1$ and $N_\alpha \subseteq N.$ Hence $N_\alpha\cap H_\kappa \in
S$ by the property of $N_\alpha.$ This is a contradiction.

This proves that $\{\alpha < \omega_1\;|\;N_\alpha \cap H_\kappa \in S\}
$ contains a club in ${\omega_1}.$

We prove now that the Projective Stationary Reflection principle implies
the Strong Reflection Principle.

Let $\kappa \geq \omega_1$ and $\theta > \kappa$. Assume that $S
\subseteq [\kappa]^\omega$ and $\theta$ is regular. Define $S^*$ to be
the following set:
 for $N\in [H_\theta]^\omega$, let $N \in S^*$ if and only if
$N\prec (H_\theta, \in, \triangle)$
and that there exists a countable $M \prec (H_\theta,\in,\triangle)$
such that $N \subseteq
M$, $N\cap \omega_1 = M\cap \omega_1$ and $M\cap \kappa \in S$ implies
that $N\cap \kappa \in S.$
\medskip
\noindent{\bf Claim}\hskip 1truecm $S^*$ is projective stationary in
$[H_\theta]^\omega.$
\medskip
Let $g:[H_\theta]^{<\omega} \to H_\theta$ and $T \subseteq \omega_1$ be
stationary in $\omega_1$.

Let $\lambda $ be a regular cardinal larger than the cardinality of
$H_\theta$.

Let $N'\prec (H_\lambda, \in ,\triangle)$ be countable such
that $N'\cap \omega_1 \in T$ and $\{\kappa, \theta, S, g\} \subseteq
N'.$ Assume that there exists a countable $M \prec (H_\theta,
\in,\triangle)$ such that $M \cap \omega_1 = N'\cap\omega_1$, $N'\cap
H_\theta \subseteq M$ and $M\cap\kappa \in S.$ (If there are no such
$M$, then $N'\cap H_\theta \in S^*$. We have what we want.) Let $N$ be
the skolem hull of $N' \cup (M\cap \kappa)$ in the structure
$(H_\lambda,\in,\triangle).$

We claim that $N\cap\kappa = M\cap \kappa.$ Hence $N\cap H_\theta \in
S^*$ and we finish the proof.

Let $\alpha \in N\cap \kappa.$ Let $\tau$ be a skolem term. Let $a \in
N'$ and $\alpha_1,\cdots, \alpha_m \in M\cap \kappa$ be such that $\alpha
= \tau(a,\alpha_1,\cdots,\alpha_m).$

Define $h : [\kappa]^m \to \kappa$ by
$$h(x_1,\cdots,x_m) = \cases{\tau(a,x_1,\cdots,x_m),& if
$\tau(a,x_1,\cdots,x_m) < \kappa,$\cr
0,&otherwise.\cr}$$
Then $h \in N'$. Hence $h \in N'\cap H_\theta \subseteq M.$ Therefore,
$\alpha = h(\alpha_1,\cdots,\alpha_m) \in M\cap \kappa.$

This finishes the proof that $S^*$ is projective stationary in
$[H_\theta]^\omega.$

Applying the Projective Stationary Reflection principle, let $\langle
N_\alpha \;|\;\alpha < \omega_1\rangle$ be an increasing continuous
$\in$--chain of length $\omega_1$ such that $N_\alpha \in S^*$ for all
$\alpha < \omega_1$. Certainly, this is what the Strong Reflection
Principle needs for the given $S$.

\hfill\square

\endsection
\st{Some Examples }

In this section, we investigate some natural examples of projective
stationary sets. These examples have been used previously by others,
and in some cases we are unsure of the authorship of these examples.

\medskip
First let us consider the saturation of the nonstationary ideal on
$\omega_1$.

Let $F$ be a set of stationary sets on $\omega_1$. $F$ is an antichain
if $A\cap B$ is nonstationary for $A \not= B$ in $F$. $F$ is a maximal
antichain if $F$ is an antichain and for every stationary subset $T
\subseteq \omega_1$ there exists some $A \in F$ such that $T \cap A$ is
stationary. The nonstationary ideal on $\omega_1$ is saturated if every
maximal antichain has size at most $\aleph_1$.

\medskip
In [7], Steel and Van Wesep showed that the nonstationary ideal may be
saturated using determinancy. In [3], Foreman, Magidor and Shelah
showed that Martin's Maximum implies that the nonstationary ideal on
$\omega_1$ is saturated. In [9], Todorcevic showed that the Strong
Reflection Principle is sufficient.
 \medskip

\example

Let $F$ be a maximal antichain of stationary subsets of $\omega_1$.
Consider the following set:
$$S = \{ N \in [H_{\omega_2}]^{\omega}\;|\;
N\ \hbox{\rm  is elementary and }\exists \,A \in F\cap N\;(N\cap
\omega_1 \in A) \}.$$

We now check that $S$ is projective stationary.

Let $T\subseteq \omega_1$ be stationary. Let $A \in F$ be such that
$A \cap T$ is stationary. Let $C$ be a club on
$[H_{\omega_2}]^{\omega}$.
Then we can find an elementary submodel $N$ of $H_{\omega_2}$ with the
property that $N\cap \omega_1 \in A\cap T$ and $N \in C$. Hence $S$ is
projective stationary.

Applying the Strong Reflection Principle, let $\{N_\alpha
\;|\;\alpha < \omega_1 \}$ be an increasing continuous
$\subseteq$--chain
of elementary submodels of $H_{\omega_2}$ such that $\alpha \subseteq
N_\alpha$ for every $\alpha$ and there exists a club $C
\subseteq \omega_1$ with that if $\alpha \in C$ then $\alpha = N_\alpha
\cap \omega_1$ and $N_\alpha \in S$.

We proceed to check that $F \subseteq X=\bigcup \{N_\alpha\;|\;\alpha <
\omega_1\}.$

Let $A \in F$. Assume that this $A \not\in X$. Let $Y$ be the skolem
hull of $X\cup \{A\}$. Let $M_\alpha$ be the skolem hull of $N_\alpha
\cup \{A\}.$ Let $D \subseteq C$ be a club such that for every $\alpha
\in D$ we have that $M_\alpha \cap \omega_1 = N_\alpha\cap \omega_1 =
\alpha.$

Let $\alpha \in D\cap A.$ Then $\alpha = N_\alpha \cap
\omega_1= M_\alpha\cap \omega_1$. There must be some $B \in F\cap
N_\alpha$ such that $\alpha \in B.$ This $B$ must be different from $A$.
Hence $A\cap B$
must be nonstationary. But by elementarity, as both $A$ and $B$ are in
$M_\alpha$, there must be some closed and unbounded subset $E \in
M_\alpha$ of $\omega_1$ such that $E \cap A \cap B$ is empty. We then
have a contradiction  because $E \in M_\alpha$ implies that
$M_\alpha\cap \omega_1 \in E$ and
$M_\alpha \cap \omega_1 \in A \cap B$.

Therefore $F \subseteq X$. Hence the size of $F$ is at most $\aleph_1.$

Hence we have the theorem of Todorcevic [1,9] that
the Strong Reflection Principle implies that
the nonstationary ideal on $\omega_1$ is saturated.

\hfill\square

\example

Our next example is to show that
assuming the Strong
Reflection Principle, every stationary set $E \subseteq \omega_2$ of
ordinals of cofinality $\omega$ contains
a closed copy of $\omega_1$. 

Fix a regular cardinal $\kappa \geq \omega_2$.

Let $E$ be a stationary subset of $\kappa$ such that every $\alpha
\in E$ has cofinality $\omega.$

Consider the following set:
$$S = \{N \in [H_{\kappa}]^{\omega}\;|\;N\ \hbox{\rm is elementary and
sup}(N\cap \kappa)\in E \}.$$

We check that $S$ is projective stationary.

Let $T \subseteq \omega_1$ be stationary. Let $f :
[H_\kappa]^{<\omega}\to H_\kappa$. Let $M \prec H_\kappa$ be such that
$\omega_1 \subseteq M$, $M$ has cardinality $\aleph_1$, $M$ is closed
under $f$ and sup$(M\cap\kappa) \in E$. Then there is a countable $N
\prec M$ such that $N$ is closed under $f$ and $N\cap\omega_1 \in T$ and
sup$(N\cap \kappa) = $ sup$(M\cap\kappa)$.

This shows that $S$ is projective stationary.

Applying the Strong Reflection Principle, let
$\{N_\alpha\;|\;\alpha < \omega_1\}$ be an increasing continuous
$\in$--chain such that
$N_\alpha \in S$ for $\alpha
\in \omega_1$. Then for each $\alpha <\omega_1$, sup$(N_\alpha\cap
\kappa) \in E.$ We are done.

Therefore we have proved the theorem of Todorcevic [1,9] that
the Strong
Reflection Principle implies that if $\kappa \geq \omega_2$
is regular, if $E \subseteq \kappa$ is stationary,
and if every $\alpha \in E$ has cofinality $\omega$, then $E$
contains a closed copy of $\omega_1$.

\hfill\square

\example

Let $E = \{\alpha < \kappa\;|\;cf(\alpha) = \omega\}.$

Let $\{E_n\;|\;n < \omega\}$ be a family of disjoint stationary subsets
of $E.$
Let $\{T_n\;|\;n < \omega\}$ be a partition of $\omega_1$ so that each
$T_n$ is stationary for $n < \omega.$

Then the following set is projective stationary:
$$S = \{N \in [H_\kappa]^{\omega}\;|\;\forall\,n< \omega\;
N\cap \omega_1 \in T_n
\Rightarrow \hbox{\rm sup}(N\cap \kappa) \in E_n \}.$$

Let $\{E_\alpha\;|\;\alpha < \omega_1\}$ be a family of disjoint
stationary subsets of $E.$
Let $\{T_\alpha\;|\;\alpha < \omega_1\}$ be a partition of $\omega_1$ so
that each $T_\alpha$ is stationary for $\alpha < \omega_1$ and for every
stationary subset $X\subseteq \omega_1$, there is some $\alpha <
\omega_1$ such that $X \cap T_\alpha$ is stationary.

Then the following set $S$ is projective stationary:
$$S = \{N \in [H_\kappa]^{\omega}\;|\;\forall\,\alpha< N\cap\omega_1\;
N\cap \omega_1 \in T_\alpha
\Rightarrow \hbox{\rm sup}(N\cap \kappa) \in E_\alpha \}.$$

To see this, let $T\subseteq \omega_1$ be stationary. Let $\alpha$ be
the least countable ordinal such that $T\cap T_\alpha$ is stationary.
Then the following set is stationary:
$$\{N\in [H_\kappa]^\omega\;|\;\alpha < N\cap \omega_1 \in
T\cap T_\alpha\,\&\,\hbox{\rm sup}(N\cap\kappa)\in E_\alpha\}.$$

This can be used to show, assuming the Strong Reflection
Principle, that if $\kappa \geq \omega_2$ is regular, then
$\kappa^{\aleph_1}=\kappa$.
By Silver's theorem [6], this in turn implies the Singular Cardinal
Hypothesis holds.
Here, by the Singular Cardinal
Hypothesis we mean the following statement:
\medskip
{\sl For every singular cardinal $\lambda$, if $2^{cf(\lambda)} <
\lambda$, then $\lambda^{cf(\lambda)} = \lambda^+.$}

Thus, we have that
the Strong Reflection Principle
implies that the Singular Cardinal Hypothesis holds.

\hfill\square

\example

We now prove that the Strong Reflection Principle implies
the Reflection Principle.
In fact, if $S \subseteq [H_\kappa]^\omega$ is
stationary then there exists an increasing continuous $\in$--chain
$\langle N_\alpha\;|\;\alpha < \omega_1\rangle$ such that $\{\alpha <
\omega_1\;|\;N_\alpha \in S\}$ is stationary in $\omega_1$.
(This is
the version of the reflection principle used by
Velickovic in [10] to show the Singular Cardinal Hypothesis.)

\medskip
Let $S \subseteq [H_\kappa]^\omega$.
There is another
projective stationary set $p(S) \subseteq [H_\kappa]^\omega$ naturally
associated with $S$.

Let $T\subseteq \omega_1$. Define $S_T=\{N \in S\;|\;N\cap \omega_1 \in
T\}.$
Then $S$ is projective stationary if and only if $S_T$ is stationary for
every stationary $T \subseteq \omega_1$.

Let $F \subseteq \{ T \subseteq \omega_1\;|\;T \hbox{\rm \ is stationary
and\ } S_T\hbox{\rm \ is not stationary }\}$ be a maximal antichain of
the smallest cardinality.

Define $p(S)$ to be the following set:
$$\{N \prec H_\kappa\;|\; N\in [H_\kappa]^\omega
\hbox{\rm \ and }
(\exists\, A \in N\cap F\,(N\cap \omega_1 \in A) \iff N \not\in S)\}.$$
Then $p(S)$ is projective stationary in $[H_\kappa]^\omega.$

If we
assume that $\langle N_\alpha\;|\;\alpha <\omega_1\rangle$ is an
increasing continuous $\in$--chain from $p(S)$, then $S$ is stationary
if and only if $\{\alpha < \omega_1\;|\;N_\alpha \in S\}$ is stationary.

Let us assume that $S$ is stationary.
Toward a contradiction, let us assume that $\{\alpha <
\omega_1\;|\;N_\alpha \in S\}$ is not stationary.
Then the complement of this set contains
a club. It follows
that for every $A \in F$
there is some $\alpha < \omega_1$ such that $A \in N_\alpha.$ Let $f :
\omega_1 \to F$ be a surjective mapping.
For each $\alpha < \omega_1$ let $C_\alpha \subseteq [H_\kappa]^\omega$
be a club to witness that $f(\alpha) \in F$. Let $C = \Delta_{\alpha<
\omega_1}C_\alpha.$ Then $C \cap S$ is stationary and the following set
$T$ is stationary in $\omega_1$:
$$T = \{N \cap \omega_1\;|\;N \in C\cap S\}.$$
Let $\alpha \in T$ be such that for all $A \in N_\alpha\cap F$ there is
a $\beta < \alpha$ with $A = f(\beta)$, for all $\beta < \alpha$,
$f(\beta) \in N_\alpha$, and $\alpha = N_\alpha \cap \omega_1$. Then
$N_\alpha \in S.$ This is a contradiction.
Hence $\{\alpha < \omega_1\;|\;N_\alpha \in S\}$ is
stationary in $\omega_1$.

To see the other direction, assume that $S$ is not stationary. By 
minimality,  the cardinality of $F$ must be one and the member of $F$
must contain a club. By elementarity, there must be such a club in
$N_0$. Therefore, for all $\alpha < \omega_1$, $N_\alpha\cap\omega_1$
must be in this club. Hence no $N_\alpha$ will be in $S$.

Furthermore, if we require the $\omega_1$ sequences to be strongly
increasing
continuous, then we can have that if $S \subseteq [H_\kappa]^\omega$ is
stationary and $T\subseteq [H_\kappa]^{\omega_1}$ is $\omega_1$--closed
and unbounded, then there exists an $X \in T$ such that $S \cap
[X]^\omega$ is stationary in $[X]^\omega$.

\hfill\square

\example

Let $C$ be a sequence defined on all the limit ordinals $< \kappa$ such
that for all limit ordinal $\alpha < \kappa$, $C_\alpha \subseteq
\alpha$ is a club in $\alpha$ and for all limit ordinals $\beta <
\alpha$, if $\beta \in C_\alpha$ is a limit point of $C_\alpha$, then
$C_\beta = \beta \cap C_\alpha.$ $C$ is a $\square^*_\kappa$ sequence if
there is no club $D \subseteq \kappa$ such that for every limit point
$\alpha$ of $D$, $C_\alpha=\alpha \cap D.$

Let $S_C$ be the set of all countable elementary submodels $N$ of
$H_\kappa$
such that, letting $\alpha = \hbox{\rm sup}(N\cap \kappa)$, $C_\alpha
\cap N$ is bounded in $\alpha$.

The following theorem is essentially due to Velickovic [10].

\theorem The following are equivalent:

(1) $S_C$ is projective stationary.

(2) $S_C$ is stationary.

(3) There is no club $D \subseteq \kappa$ such that for all limit point
$\alpha$ of $D$, $C_\alpha = \alpha \cap D.$ That is, $C$ is a
$\square^*_\kappa$ sequence.

\proof

(2) $\Rightarrow$ (3).
Assume that (3) is false. Let $D$ be a witness. Let $\theta > \kappa$ be
a sufficiently large regular cardinal. Let $N\prec H_\theta$ be a
countable elementary submodel containing all the relevant objects. Then
$N\cap H_\kappa \not\in S_C.$

(3) $\Rightarrow$ (1).  This is essentially due to Velickovic [10].

For an $f : [H_\kappa]^{<\omega} \to H_\kappa$, and a countable
ordinal $\alpha$, let us consider the following game $G(\alpha, f)$, due
to Velickovic [10].
Two players, Player One and Player Two, take in turn playing in $\omega$
many steps. Player One plays an interval $I_n \subseteq \kappa$ of size
less than $\kappa$ and an ordinal $\alpha _n \in I_n$ at his $n$th move.
Player Two plays an ordinal $\beta_n$ in his responding move. At the end
of the play, Player One wins the play if and only if for every $n <
\omega$, $\beta_n < \hbox{\rm inf}(I_{n+1})$ and the elementary submodel
$N$ of $H_\kappa$,
which is generated by $\alpha \cup\{\alpha_n\;|\;n < \omega\}$ and which
is also closed under $f$,
has the property that $N \cap \omega_1 = \alpha$ and $N \cap \kappa
\subseteq \bigcup_{n<\omega} I_n.$

Notice that if Player One loses a play, Player One loses at some stage
even before the game is over. Hence this game is a determined game.

\Lemma{(Velickovic)} There exists a club $E_f\subseteq \omega_1$ such
that for every $\alpha \in E_f$, Player One has a winning strategy
$\sigma_f^\alpha$ for the game $G(\alpha,f).$
\vskip 5truept
Let $f : [H_\kappa]^{<\omega}\to H_\kappa$. Let $T\subseteq \omega_1$ be
stationary.

Assume (3), we need to find a countable elementary submodel $N\prec
H_\kappa$ with the property that $N$ is closed under $f$, $N\cap
\omega_1 \in T$ and $N \in S_C.$

Let $\alpha \in T \cap E_f$. Let $\sigma$ be a winning strategy for
Player One in the game $G(\alpha, f)$.

Let $\theta = (2^\kappa)^+.$ Let $\langle M_\alpha\;|\;\alpha <
\kappa\rangle$ be an increasing continuous $\in$--chain of elementary
submodels of $H_\theta$ of size $< \kappa$ such that all the relevant
objects are in $M_0$ and $M_\alpha \cap \kappa \in \kappa$ for all
$\alpha
< \kappa.$ Let $D = \{M_\alpha \cap \kappa \;|\;\alpha < \kappa\}.$
Then $D$ is a club in $\kappa.$

Let $M \prec H_\theta$ be an elementary submodel of size less than
$\kappa$, containing all the relevant objects, such that $M\cap \kappa$
is an ordinal of cofinality $\omega_1.$ Let $\delta = M\cap\kappa.$

Assume that there is a $X \in M$ such that $\delta \cap X \subseteq
C_\delta$ and $M \models `` X$ is a club in $\kappa$ ''. Let $Y$ be the
set of all the limit points of $X$. Then $Y \in M$. Let $\gamma < \beta$
be in $M$ such that both are in $Y$. Then $\gamma$ and $\beta$ are limit
points of $C_\delta$. Hence $C_\gamma = \gamma \cap C_\beta.$ Therefore,
for every pair $\gamma < \beta$ from $Y$, we have $C_\gamma = \gamma
\cap C_\beta.$ This contradicts to our assumption (3).

Thus, for every club $X \subseteq \kappa$, if $X \in M$ then there is
some ordinal $\gamma \in X\cap \delta$ which is not in $C_\delta.$
It follows that for every $\gamma < \delta$, there is some $\xi <
\delta$ such that $\eta_\xi = M_\xi\cap \kappa$ is not in $C_\delta$ and
$\eta_\xi > \gamma.$

Inductively pick $\xi_n$ so that $\eta_{\xi_n} \not\in C_\delta$ and
there is some ordinal $\beta \in C_\delta$ with $\eta_{\xi_n} < \beta <
\eta_{\xi_{n+1}}.$

Let $\beta = \bigcup\{\eta_{\xi_n}\;|\;n < \omega\}$. Then $C_\beta =
\beta \cap C_\delta.$

Consider the
following play of the game $G(\alpha, f)$.
Player One follows his winning strategy $\sigma$ which is in $M_0$.
Player Two plays $\mu_n = $ max$(\eta_{\xi_n}\cap C_\delta)$, which is
in $M_{\xi_n}$, in his $n$th move.
Let $I_n$ and $\beta_n$ be the $n$th move of Player One following his
winning strategy $\sigma$. Then $\beta_n \in I_n \subseteq (\mu_n,
\eta_{\xi_n})$ for all $n < \omega$.
Let $N\prec M$ be generated by $\alpha \cup \{ \beta_n\;|\;n<\omega\}$
which is closed under $f$. Then $N \cap \omega_1 \in T$ and sup$N\cap
\kappa = \beta$ and $N \cap C_\beta \subseteq I_0$. Hence $N \cap
H_\kappa \in S_C.$

This shows that (3) implies (1).

\hfill\square

\corollary If for every stationary $S\subseteq [H_\kappa]^\omega$ there
is an increasing continuous $\in$--chain
$\langle N_\alpha\;|\;\alpha < \omega_1\rangle$ such that $\{\alpha
<\omega_1\;|\;N_\alpha \in S\}$ is stationary, then there is a club
$C\subseteq \kappa$ such that for every limit point $\alpha \in C$,
$C_\alpha = \alpha \cap C.$ In particular, there is no
$\square^*_\kappa$ sequence.

\proof
By the previous theorem, we need only to prove that $S_C$ is not
stationary.

Toward a contradiction, let us assume that $S_C$ is stationary.
Let $\langle N_\alpha\;|\;\alpha < \omega_1\rangle$ be an increasing
continuous $\in$--chain be such that $T = \{\alpha <
\omega_1\;|\;N_\alpha \in S_C\}$ is stationary.
Let $\gamma_\alpha = \hbox{\rm \ sup}(N_\alpha\cap\kappa)$ for
$\alpha<\omega_1$. Let $\delta = \hbox{\rm \
sup}\{\gamma_\alpha\;|\;\alpha < \omega_1\}.$
Let $\alpha \in T$ be such that $\gamma_\alpha \in C_\delta$ and
$\gamma_\alpha$ is a limit point of
$C_\delta\cap\{\gamma_\beta\;|\;\beta < \omega_1\}.$ We get a
contradiction, since $\alpha \in T$ implies that $N_\alpha \in S_C$ and
$C_{\gamma_\alpha}=\gamma_\alpha\cap C_\delta$ whose intersection with
$N_\alpha$ is not bounded in $\gamma_\alpha.$

\hfill\square

\endsection

\st {Saturation and Reflection}

In this section, we take a closer look at the saturation of the
nonstationary ideal on $\omega_1$. It turns out that it itself is a kind of
reflection property.

\def\filter{\cal F}
\def\Rarrow{\Rightarrow}

To start with, let us define first certain filters. For a
regular cardinal $\kappa \geq \omega_2$, for $X \subseteq
[H_\kappa]^\omega$, let $X$ be in $\filter_\kappa$ if and only if for
every stationary subset $A \subseteq \omega_1$ there exist a stationary
subset $B \subseteq A$ and a closed and unbounded subset $C \subseteq
[H_\kappa]^\omega$ such that $\{N \in C \;|\;N\cap \omega_1 \in B \}
\subseteq X.$

\theorem The following are equivalent:

(1) The nonstationary ideal $NS_{\omega_1}$ on $\omega_1$ is saturated.

(2) For every regular cardinal $\kappa \geq \omega_2$, for every
stationary set $S \subseteq [H_\kappa]^\omega$, there exists a
stationary set $A \subseteq \omega_1$ such that $S$ is $A$--projective
stationary (i.e., for every stationary $B \subseteq A$, the set $\{N \in
S \;|\;N\cap \omega_1 \in B\}$ is stationary).

(3) For every regular cardinal $\kappa \geq \omega_2$, for every $X
\subseteq [H_\kappa]^\omega$, $X \in \filter_\kappa$ if and only if $X$
contains a closed and unbounded subset. Namely, the filter
$\filter_\kappa$ is just the club filter on $[H_\kappa]^\omega$.

(4) For every regular cardinal $\kappa \geq \omega_2$, for every $X \in
\filter_\kappa$, there exists an increasing continuous $\in$--chain
$\langle N_\alpha \;|\;\alpha < \omega_1\rangle$ of countable
elementary submodels of $H_\kappa$ of length
$\omega_1$ such that $N_\alpha \in X$ for all $\alpha < \omega_1$.

\proof  (1) $\Rarrow$ (2)

Let $S \subseteq [H_\kappa]^\omega$ be a stationary set. If $S$ is
projective stationary, then there is nothing needed to be proved. So we
assume that $S$ is not projective stationary.

For a stationary $T \subseteq \omega_1$, we let $T \in F$ if and only
if there exists a club $C \subseteq [H_\kappa]^\omega$ such that
$$T \cap \{\,N\cap\omega_1\;|\;N \in C\cap S\} = \emptyset.$$
Since $S$ is not projective stationary, $F$ is not empty.

Because the nonstationary ideal on $\omega_1$ is saturated, we can take
a sequence $\{A_\alpha\;|\;\alpha < \omega_1\}$ from $F$ so that for
every $A \in F$ , $A - \nabla_{\alpha < \omega_1}A_\alpha$ is
nonstationary, where $\nabla_{\alpha< \omega_1}A_\alpha$ is the diagonal
union of the sequence $\{A_\alpha\;|\;\alpha<\omega_1\}$ defined by
$$\nabla_{\alpha< \omega_1}A_\alpha = \{\beta<\omega\;|\;\exists\,\alpha
< \beta\;(\beta \in A_\alpha)\}.$$

\def\dua{\nabla_{\alpha< \omega_1}A_\alpha}
\def\dic{\Delta_{\alpha<\omega_1}C_\alpha}

For each $\alpha <\omega$, we let $C_\alpha \subseteq [H_\kappa]^\omega$
be a witness for $A_\alpha$ to be in $F$.

Let $\dic$ be the diagonal intersection of the sequence $\{C_\alpha\;|\;
\alpha < \omega_1\}$ defined by
$$\dic = \{N \in [H_\kappa]^\omega\;|\;\forall\,\alpha \in
N\cap\omega_1\; N \in C_\alpha\}.$$

Let $A = \omega - \dua$. We claim that $S$ is $A$--projective
stationary.

First notice that $A$ is stationary. This is because $\dic \cap S$ is
stationary.
Now let $T \subseteq A$ be stationary and let
$C \subseteq [H_\kappa]^\omega$ be a club.

Since $T \cap \dua$ is not stationary,
$T \not\in F$. Therefore the following intersection is not empty:
$$T\cap\{N\cap\omega_1\;|\;N \in C\cap\dic \cap S \}.$$
Hence $\{ N \in S \;|\; N\cap \omega_1 \in T \}$ is stationary.

\medskip
(2) $\Rarrow$ (3)

Since every closed and unbounded subset of $[H_\kappa]^\omega$ is in the
filter $\filter_\kappa$, and (2) implies that every stationary subset of
$[H_\kappa]^\omega$ intersects every member of the filter
$\filter_\kappa$, the filter $\filter_\kappa$ and the club filter on
$[H_\kappa]^\omega$ are the same filter.
\medskip
(4) $\Rarrow$ (1)

Let $F$ be a maximal antichain of stationary subsets of $\omega_1$.
Let $S_F = \{ N \in [H_\kappa]^\omega\;|\;\exists\;A \in F \cap
N\;(N\cap\omega_1 \in A) \}.$ Then $S_F \in \filter_\kappa$. Namely, for
a given stationary subset $A \subseteq \omega_1$, let $T \in F$ be such
that $B = A \cap T $ is stationary. Let $C$ be the club of all countable
elementary submodels of $H_\kappa$ which contains $T$. Then if $N \in C$
and $N\cap\omega_1 \in B$ then $N \in S_F$. Now by the proof of theorem
2.1, assuming (4), $F$ has cardinality at most $\aleph_1$.

This completes the proof.

\hfill\square
\bigskip
The next thing we want to show is that the presaturation of the
nonstationary ideal on $\omega_1$ is also a kind of reflection property,
which in turn is equivalent to the filter $\filter_\kappa$ being
countably closed.

First let us recall that the nonstaionary ideal on $\omega_1$ is {\bf
presaturated} if for every countable sequence $\{F_n \;|\;n<\omega\}$ of
maximal antichains of the nonstationary ideal on $\omega_1$, for every
stationary subset $T$, there exists a stationary subset $B \subseteq T$
such that for each $n$ the set $\{A \in F_n\;|\;A\cap B$ is stationary
$\}$ has cardinality at most $\aleph_1$ (see [3]).

\theorem The following are equivalent:

(1) The nonstationary ideal $NS_{\omega_1}$ on $\omega_1$ is
presaturated.

(2) For every regular $\kappa \geq \omega_2$, the filter
$\filter_\kappa$ is $\sigma$--closed.

(3) For every regular cardinal $\kappa \geq \omega_2$, for every
countable sequence $\langle X_n \;|\; n < \omega\rangle$ from
$\filter_\kappa$, for every stationary subset $T \subseteq \omega_1$,
there exists an $M \prec H_\kappa$ such that $\omega_1 \subseteq M$, $M$
has cardinality $\aleph_1$ and $\{N \in [M]^\omega \cap
\bigcap_{n<\omega}X_n\;|\; N \cap\omega_1 \in T\}$ is stationary in
$[M]^\omega.$

\proof
(1) $\Rarrow$ (2)

Let $X_n \in \filter_\kappa$ for $n < \omega$. Let $X =
\bigcap_{n<\omega}X_n.$ We need to show that $X \in \filter_\kappa.$

For each $n < \omega$, let $\langle B^n_\alpha,\,C^n_\alpha\;|\;\alpha <
\theta_n\rangle$ be such that all the $B^n_\alpha$'s form a maximal
antichain of stationary sets and for each $\alpha < \theta_n$, $N \in
C^n_\alpha$ and $N\cap\omega_1 \in B^n_\alpha$ imply that $N \in X_n$.

Let $A\subseteq \omega_1$ be stationary. Applying the presaturation
property of the nonstationary ideal, let $B \subseteq A$ be stationary
such that for each $n < \omega$, the set $I_n = \{\alpha < \theta_n
\;|\; B \cap B^n_\alpha \hbox{\rm\ is stationary }\}$ has cardinality at
most $\aleph_1$. To simplify the notation, we assume that $I_n =
\omega_1$. Let $$T = B \cap \bigcap_{n<\omega}\nabla_{\alpha<
\omega_1}B^n_\alpha$$ and let
$$C = \bigcap_{n < \omega}{\Delta_{\alpha<\omega_1}C^n_\alpha}.$$
It follows that if $N \in C$ and $N\cap\omega_1 \in T$ then $N \in X.$

Hence the filter $\filter_\kappa$ is $\sigma$--closed.
\medskip
(2) $\Rarrow$ (3)

It will be sufficient to show that every $X \in \filter_\kappa$
reflects.

Let $X \in \filter_\kappa$. Let $T$ be a stationary subset of
$\omega_1$. Let $B \subseteq T$ be stationary
and let $C \subseteq [H_\kappa]^\omega$ be a closed and unbounded subset
such that for all $N \in C$, $N\cap \omega_1\in B$ implies that $N \in
X$.

Let $\langle N_\alpha \;|\;\alpha < \omega_1\rangle$ be an increasing
continuous $\in$--chain from $C$. Let $M = \bigcup_{\alpha <
\omega_1}N_\alpha$. Then $\{N_\alpha\;|\;\alpha = N_\alpha\cap \omega_1
\in B\}$ is stationary in $[M]^\omega$.

\medskip
(3) $\Rarrow$ (1)

Fix a sequence $\{ F_n\;|\;n < \omega\}$ of maximal antichains of the
nonstationary ideal on $\omega_1$.
For each $n < \omega$, let
$$S_n = \{N\prec H_{\kappa}\;|\;N\ \hbox{\rm is countable and\ }
\exists\,A \in F_n\cap N\;(N\cap \omega_1 \in A)\}.$$
Then every $S_n\in \filter_\kappa$.

Fix a stationary subset $T$ of $\omega_1$.
Applying (3), let $X \prec H_{\omega_2}$ be such that
$\omega_1\subseteq X$, $X$ has cardinality $\aleph_1$ and
$\{N \in [X]^\omega\cap \bigcap_{n<\omega}S_n\;|\;N\cap\omega_1 \in T\}$
is stationary in $[X]^\omega.$
Write $X = \bigcup_{\alpha<\omega_1} N_\alpha$, the union of an
increasing
continuous $\subseteq$--chain of countable elementary submodels with
$\alpha \subseteq N_\alpha$. Then
$$T_0 = \{ \alpha \in T\;|\;\alpha = N_\alpha \cap
\omega_1\;\&\;N_\alpha
\in \bigcap_{n<\omega}S_n\}$$ is a stationary subset of $T$.

We claim that this $T_0$ has the desired property.

Assume not. Let $n$ be the least counterexample. Let $A \in F_n - X$ be
such that $A \cap T_0$ is stationary. Let $M_\alpha $ be the skolem hull
of $N_\alpha \cup \{A\}$ for each $\alpha < \omega_1$. Then
$$C=\{\alpha < \omega_1\;|\; \alpha =
N_\alpha\cap\omega_1=M_\alpha\cap\omega_1\}$$
is a club. Let $\alpha \in C\cap A\cap T_0.$ Let $B \in N_\alpha \cap
F_n$ be such that $\alpha \in B.$ Since both $A$ and $B$ are in
$M_\alpha$ which is an elementary submodel of $H_{\kappa}$, and
$M_\alpha\cap\omega_1 = \alpha \in A\cap B$, we conclude that $A\cap B$
is stationary. This is a contradiction.

Hence for every $n < \omega$, for every $A \in F_n$, if $A\cap T_0$ is
stationary, then $A \in X$. We are done since $X$ has cardinality
$\aleph_1$.
\hfill\square
\endsection

\vfill\eject

\reference

\referno=1

\ref{M. Bekkali}{}{Topics in Set Theory}
{Lecture Notes in Mathematics, Vol. 1476, Springer--Verlag, Berlin, New
York, 1991 }

\ref{Q. Feng and T. Jech}{Local Clubs, Reflection, and Preserving
Stationary Sets}{Proc. London Math. Soc.}{(3) 58 (1989) 237--257}

\ref{M. Foreman, M. Magidor and S. Shelah}{Martin's Maximum, Saturated
Ideals, and Nonregular Ultrafilters. Part I}{Annals of Mathematics}{127
(1988) 1--47}

\ref{T. Jech}{}{Set Theory}{Academic Press, New York 1978}

\ref{T. Jech}{}{Multiple Forcing}{Cambridge Tracts in Mathematics,
Cambridge University Press, 1986}

\ref{J. Silver}{On the Singular Cardinals Problem}{Proc. Internat'l
Cong. Math.}{Vancouver, B. C., 1974, Vol. 1, 265--268}

\ref{J. Steel and R. van Wesep}{Two Consequences of Determinancy
Consistent with Choice}{Trans. A. M. S.}{272(1) (1982) 67--85}

\ref{S. Todorcevic}{Reflecting Stationary Sets I}{}{Hand Written Notes,
1985}

\ref{S. Todorcevic}{Strong Reflections}{}{Hand Written Notes, 1987}

\ref{B. Velickovic}{Forcing Axioms and Stationary Sets}{Advances in
Mathematics}{94 (1992) 256--284}

\vskip 20truept
\noindent Department of Mathematics, National University of Singapore,
Singapore 0511

\smallskip
\noindent{\sl Email:} {\tt matqfeng@leonis.nus.sg}
\medskip
\noindent Department of Mathematics, Pennsylvania State University,
University Park, PA 16802, USA
\smallskip
\noindent{\sl Email:} {\tt jech@math.psu.edu}

\vfill\eject

\pageno=1
\input vanilla.sty

\title
{{\bf Projective Stationary Sets and Strong Reflection Principle}}
\endtitle
\par
\author
Qi Feng and Thomas Jech
\endauthor
\par
\heading
Abstract
\footnote{AMS classification: 03E55, 03E65}
\footnote"*"{Keywords: Stationary Sets, Reflection}
\endheading

{We study projective stationary sets.
The Projective Stationary Reflection principle is  the statement that
every projective stationary set contains an increasing
continuous $\in$--chain of length $\omega_1$.
We show that if Martin's Maximum
holds, then the Projective Stationary Reflection Principle holds. Also
it is equivalent to the Strong Reflection Principle. We show that
the saturation of the nonstationary ideal on $\omega_1$ is
equivalent to a certain kind of reflection.}

\bye